\documentclass[12pt,a4paper]{article}
\usepackage{amsmath}
\usepackage{graphicx}
\allowdisplaybreaks[1]

\newcommand{\pr}{\rightarrow}

\newcommand{\ba}{\begin{array}}
\newcommand{\ea}{\end{array}}

\newcommand{\vart}{\vartheta}
\newcommand{\varp}{\varphi}
\newcommand{\eps}{\varepsilon}

\newenvironment{inspring}[1]%
{\begin{list}{}{\setlength{\rightmargin}{0cm}
                \setlength{\listparindent}{0cm}
                \settowidth{\labelwidth}{\mbox{#1}}
                \setlength{\leftmargin}{1.1\labelwidth}
                \setlength{\labelsep}{.1\labelwidth}}}%
{\end{list}}

\newcommand{\bi}[1]{\begin{inspring}{#1}}
\newcommand{\ei}{\end{inspring}}

\newcommand{\beq}{\begin{equation}}
\newcommand{\eq}{\end{equation}}
\catcode`\@=11

\font\tenmsa=msam10 \font\sevenmsa=msam7 \font\fivemsa=msam5
\font\tenmsb=msbm10 \font\sevenmsb=msbm7 \font\fivemsb=msbm5
\newfam\msafam
\newfam\msbfam
\textfont\msafam=\tenmsa  \scriptfont\msafam=\sevenmsa
  \scriptscriptfont\msafam=\fivemsa
\textfont\msbfam=\tenmsb  \scriptfont\msbfam=\sevenmsb
  \scriptscriptfont\msbfam=\fivemsb

\def\Bbb{\ifmmode\let\next\Bbb@\else
 \def\next{\errmessage{Use \string\Bbb\space only in math mode}}\fi\next}
\def\Bbb@#1{{\Bbb@@{#1}}}
\def\Bbb@@#1{\fam\msbfam#1}
\newcommand{\dR}{{\Bbb R}}

\numberwithin{thm}{section}

\title{Sharp bound on the radial derivatives of the Zernike circle polynomials (disk polynomials)}
\author{A.J.E.M.\ Janssen \\
Eindhoven University of Technology, \\
Department of Mathematics and Computer Science, \\
P.O.\ Box 513, 5600 MB Eindhoven, The Netherlands. \\
E-mail address: \{a.j.e.m.janssen@tue.nl\} \\
Tel.: +31-402474541}
\date{}

\begin{document}
\maketitle
\mbox{} \\ \\ \\ \\ \\
\noindent
{\bf Abstract.} \\
We sharpen the bound $n^{2k}$ on the maximum modulus of the $k^{{\rm th}}$ radial derivative of the Zernike circle polynomials (disk polynomials) of degree $n$ to $n^2(n^2-1^2)\cdot ... \cdot(n^2-(k-1)^2)/2^k(1/2)_k$. This bound is obtained from a result of Koornwinder on the non-negativity of connection coefficients of the radial parts of the circle polynomials when expanded into a series of Chebyshev polynomials of the first kind. The new bound is shown to be sharp for, for instance, Zernike circle polynomials of degree $n$ and azimuthal order $m$ when $m=O(\sqrt{n})$ by using an explicit expression for the connection coefficients in terms of squares of Jacobi polynomials evaluated at 0. 
\newpage\noindent
\section{Introduction} \label{sec1}
\mbox{} \\[-9mm]

The Zernike circle polynomials (ZCPs) $Z_n^m$ are defined for integer $n$ and $m$ such that $n-|m|$ is even and non-negative by
\beq \label{e1}
Z_n^m(r\,e^{i\varp})=e^{im\varp}\,R_n^{|m|}(r)=e^{im\varp}\,r^m\,P^{(0,|m|)}_{\frac{n-|m|}{2}}(2r^2-1)~,
\eq
where $0\leq r\leq 1$, $\varp\in\dR$, and where $R_n^{|m|}(r)$ is its radial part with $P_p^{(\alpha,\beta)}(x)$ the Jacobi polynomial of degree $p=0,1,...$ corresponding to the weight function $(1-x)^{\alpha}(1+x)^{\beta}$, $-1\leq x\leq1$. The Zernike circle polynomials are directly related to the disk polynomials $R_{k,l}^{(\alpha)}$, see \cite{ref8}, 18.37(i) on pp.~477--478, by $Z_n^m=R_{\frac12(n+m),\frac12(n-m)}^{(0)}$.

The ZCPs were introduced by Zernike \cite{ref10} in connection with his phase contrast method and further elaborated by Nijboer in his thesis \cite{ref7} on the diffraction theory of aberrations for optical systems with circular pupils. The ZCPs were also used by Cormack \cite{ref2} in the context of computerized tomography with functions on the disk to be reconstructed from their line integrals. The usefulness of the ZCPs in these fields arises from the fact that they are a complete orthogonal set of functions on the unit disk that have a closed-form expression for both their 2-D Fourier transform and their Radon transform. For the basic properties of ZCPs and their application in optical diffraction theory, see \cite{ref1}, Appendix~VII and Chapter~9, Section~2. For further mathematical properties of the ZCPs, among which an addition-type formula, see \cite{ref3}.

In a recent statistical study of the Radon transform in a medical imaging context, there is interest in upper bounds for the modulus of the radial derivatives of the ZCPs, see \cite{ref6}, Appendix~B, where it is shown that for $k=0,1,...$
\beq \label{e2}
\max_{0\leq r\leq1}\,|(R_n^{|m|})^{(k)}(r)|\leq n^{2k}~,~~~~~~0\leq r\leq1~,
\eq
with $f^{(k)}(r)$ denoting the $k^{{\rm th}}$ derivative of $f$ at $r$. We shall sharpen this bound to
\beq \label{e3}
\max_{0\leq r\leq1}\,|(R_n^{|m|})^{(k)}(r)|\leq\frac{n^2(n^2-1^2)\cdot ...\cdot(n^2-(k-1)^2)}{2^k(1/2)_k}~,~~~~~~0\leq r\leq1~,
\eq
where $(\alpha)_0=1$, $(\alpha)_k=\alpha(\alpha+1)\cdot...\cdot(\alpha+k-1)$, $k=1,2,...$ (Pochhammer symbol), and where the right-hand side of (\ref{e3}) is to be interpreted as 1 for the case that $k=0$. Observe that the right-hand side of (\ref{e3}) vanishes when $k>n$. The inequality $|R_n^{|m|}(r)|\leq1$, $0\leq r\leq1$, i.e., the case $k=0$ in (\ref{e3}), follows from non-negativity of the connection coefficients that occur when $R_n^{|m|}(r)$ is expanded in a series of Chebyshev polynomials of the first kind, as was shown by Koornwinder \cite{ref5}, and the fact that $R_n^{|m|}(1)=1$. The general case $k=1,2,...$ follows by explicit results for derivatives of Chebyshev polynomials in terms of Gegenbauer polynomials, having an explicit expression for the maximum modulus of these. The inequality in (\ref{e3}) is sharp, in the sense that for large $n$ and $m=O(\sqrt{n})$, the inequality is obtained within a factor of the order $1/\sqrt{n}$. This is shown by using a result of Janssen \cite{ref3}, Section~5, that gives the connection coefficients when expanding $R_n^{|m|}$ into a series of Chebyshev polynomials in an explicit form.

\section{Preliminary result} \label{sec2}
\mbox{} \\[-9mm]

In \cite{ref3}, Section~5, it is shown that
\beq \label{e4}
R_n^{|m|}(\cos\vart)=\sum_{j=0}^{\lfloor n/2\rfloor}\,a_{n-2j}\cos(n-2j)\,\vart~,~~~~~~0\leq\vart\leq\frac{\pi}{2}~,
\eq
where $\lfloor n/2\rfloor$ is the largest integer not exceeding $n/2$, and, for non-negative integer $i$ with $n-i$ non-negative and even,
\beq \label{e5}
a_i=a_i(m)=\eps_i\,\frac{p!\,q!}{s!\,t!}\,(\tfrac12)^l\,(P_p^{(\gamma,\eps)}(0))^2~,
\eq
in which $\eps_0=1$, $\eps_1=\eps_2=\cdots =2$, and
\beq \label{e6}
p=\frac{n-l}{2}\,,~~q=\frac{n+l}{2}\,,~~s=\frac{n-r}{2}\,,~~t=\frac{n+r}{2}\,,~~\gamma=\frac{l-r}{2}\,,~~ \delta=\frac{l+r}{2}~,
\eq
with
\beq \label{e7}
l=\max(|m|,i)~,~~~~~~r=\min(|m|,i)~.
\eq
Note that $p$, $q$, $s$, $t$, $\gamma$, $\delta$ in (\ref{e6}) are non-negative and integer since $n$, $m$, $i$ have same parity. The $a_i$ are non-negative, and this confirms Koornwinder's result. Since $R_n^{|m|}(1)=P_{\frac12(n-|m|)}^{(0,|m|)}(1)=1$, we have
\beq \label{e8}
|R_n^{|m|}(\cos\vart)|\leq\sum_{j=0}^{\lfloor n/2\rfloor}\,a_{n-2j}=R_n^{|m|}(1)=1~,
\eq
showing (\ref{e3}) for $k=0$.

The result (\ref{e4}--\ref{e7}) gives specific information as to when $a_i=0$. For instance, when $i=0$, we have $l=|m|$, $r=0$, and so $\gamma=\delta$, $p=\tfrac12(n-|m|)$. Therefore $a_0=0\Leftrightarrow \tfrac12(n-|m|)$ is odd, since $P_p^{(\gamma,\gamma)}(0)=0\Leftrightarrow p$ is odd, see \cite{ref8}, second line of Table~18.61 on p.~444. Similarly, when $m=0$ we have $a_i=0$ if and only if $\tfrac12(n-i)$ is odd.

An important case of a non-vanishing $a_{n-2j}$ is $a_n=a_n(m)$ (so $j=0$) that is used in Section~\ref{sec3} where we discuss sharpness of the bound in (\ref{e3}). Then $l=i=n$, $r=|m|$, $p=0$, and since $P_0^{(\gamma,\delta)}(0)=1$ , we thus get
\beq \label{e9}
a_n=2(\tfrac12)^n\,\Bigl(\!\!\ba{c}
n \\[2mm] \dfrac{n-|m|}{2}
\ea\!\!\Bigr)~,~~~~~~n>0~.
\eq
In a similar manner, see \cite{ref9}, bottom of p.~161 for $P_1^{(\gamma,\delta)}(0)$ and $P_2^{(\gamma,\delta)}(0)$,
\begin{eqnarray} \label{e10}
& \mbox{} & a_{n-2}=a_n\,\dfrac{|m|^2}{n}~,~~~~~~n-2>0\,,~~n-2\geq|m|~; \nonumber \\[2mm]
& & a_{n-4}=a_n\,\dfrac{(n-|m|^2)^2}{2n(n-1)}~,~~~~~~n-4>0\,,~~n-4\geq|m|~.
\end{eqnarray}

\section{Proof of the new bound} \label{sec3}
\mbox{} \\[-9mm]

Let $k=1,2,...\,$. Setting $r=\cos\vart$ in (\ref{e4}) and using that $\cos(n-2j)\,\vart=T_{n-2j}(\cos\vart)$, with $T_l$ the Chebyshev polynomial of the first kind and of degree $l=0,1,..\,$, we have
\beq \label{e11}
(R_n^{|m|})^{(k)}(r)=\sum_{j=0}^{\lfloor n/2\rfloor}\,a_{n-2j}\,T_{n-2j}^{(k)}(r)~.
\eq
From \cite{ref8}, 18.9.21 on p.~447, 18.7.4 on p.~444 and 18.9.19 on p.~446 (used $k-1$ times), we have
\beq \label{e12}
T_{n-2j}^{(k)}(r)=2^{k-1}(k-1)!\,(n-2j)\,C_{n-2j-k}^k(r)~,
\eq
where $C_i^k$ is the Gegenbauer polynomial of degree $i$ corresponding to the weight function $(1-x^2)^{k-1/2}$, $-1\leq x\leq1$ with $C_i^k=0$ for $i={-}1,{-}2,...\,$. The result (\ref{e12}) can be found in this form in \cite{ref9}, first item in line 6 of p.~188, from which we borrow the notation $C_i^k$, rather than the notation $C_i^{(k)}$ as used in \cite{ref8}, Chapter~18, to avoid confusion with the operation of taking the $k^{{\rm th}}$ derivative. From \cite{ref8}, 18.14.4 on p.~450, we have that $|C_i^k(r)|$ is maximal (strictly when $i>0$) at $r={\pm}1$, with maximum value
\beq \label{e13}
C_i^k(1)=\frac{(2k)_i}{i!}~.
\eq
Since $a_{n-2j}\geq0$, we have from (\ref{e11}) and (\ref{e12}) that $|(R_n^{|m|})^{(k)}(r)|$ is maximal at $r=1$, with maximal value
\begin{eqnarray} \label{e14}
(R_n^{|m|})^{(k)}(1) & = & \sum_{j=0}^{\lfloor n/2\rfloor}\,a_{n-2j}\,T_{n-2j}^{(k)}(1) \nonumber \\[3.5mm]
& = & 2^{k-1}(k-1)!\,\sum_{j=0}^{\lfloor n/2\rfloor}\,a_{n-2j}(n-2j)\,\frac{(2k)_{n-2j-k}}{(n-2j-k)!}~.
\end{eqnarray}
The maximum at $r=1$ is strict when $n>0$ since, see (\ref{e9}), $a_n>0$.

The quantity $i(2k)_{i-k}/(i-k)!$ is positive and increasing in $i=k,k+1,...\,$, while $\sum_{j=0}^{\lfloor n/2\rfloor}\,a_{n-2j}=1$ and $a_{n-2j}\geq0$ for $j=0,1,...,\lfloor n/2\rfloor$. Therefore
\beq \label{e15}
2^{k-1}(k-1)!\,n\,\frac{(2k)_{n-k}}{(n-k)!}\,a_n\leq (R_n^{|m|})^{(k)}(1)\leq 2^{k-1}(k-1)!\,n\,\frac{(2k)_{n-k}}{(n-k)!}~.
\eq
The third member in (\ref{e15}) can be rewritten as
\begin{eqnarray} \label{e16}
2^{k-1}(k-1)!\,n\,\frac{(2k)_{n-k}}{(n-k)!} & = & \frac{2^{k-1}(k-1)!}{(2k-1)!}\,n\,\frac{(n+k-1)!}{(n-k)!} \nonumber \\[3mm]
& = & \frac{n^2(n^2-1^2)\cdot ...\cdot(n^2-(k-1)^2)}{2^k(1/2)_k}~,
\end{eqnarray}
and this gives the upper bound (\ref{e3}) as well as a lower bound for \\
$\max_{0\leq r\leq1}\, |(R_n^{|m|})^{(k)}(r)|$.

\section{Sharpness of the bounds in (\ref{e15})} \label{sec4}
\mbox{} \\[-9mm]

We have
\beq \label{e17}
(R_n^{|m|})^{(k)}(1)=\sum_{j=0}^{\lfloor n/2\rfloor}\,a_{n-2j}(m)\,T_{n-2j}^{(k)}(1)~,
\eq
where $T_1^{(k)}(1)$ are given for $i=k,k+1,...$ by
\beq \label{e18}
T_i^{(k)}(1)=\frac{i^2(i^2-1)\cdot ...\cdot(i^2-(k-1)^2)}{2^k(1/2)_k}= 2^{k-1}(k-1)!\, \frac{i(2k)_{i-k}}{(i-k)!}~,
\eq
and vanish for $i<k$. When (for a given $n$ and $m$) $k$ increases, the terms $a_{n-2j}(m)$ in the series in (\ref{e7}) with small $j$ get relatively more weight than those with large $j$. It is therefore expected that the inquality in (\ref{e3}), i.e., the upper bound in (\ref{e15}), tends to be sharp for small $k$ while the lower bound in (\ref{e15}) tends to be sharp for large $k$. Indeed, there is equality in the second inquality in (\ref{e15}) when $k=0$ and in the first inequality in (\ref{e15}) when $k=n,n-1$ (for then the only non-vanishing term in the series in (\ref{e17}) is the one with $j=0$).

For the cases $k=1,2\,$, we have
\beq \label{e19}
(R_n^{|m|})'(1)=|m|+2w~,~~~~~~(R_n^{|m|})''(1)=|m|(|m|-1)+2w(w+|m|-1)~,
\eq
where $w=\frac14(n-|m|)(n+|m|+2)$. With $n^2$, $\frac13\,n^2(n^2-1)$ being the right-hand side of (\ref{e3}) for $k=1,2\,$, we have from (\ref{e19})
\beq \label{e20}
\frac{(R_n^{|m|})'(1)}{n^2}\pr \frac12~,~~~~~~\frac{(R_n^{|m|})''(1)}{\tfrac13\,n^2(n^2-1)}\pr\frac38
\eq
when $n\pr\infty$ and $m=o(n)$. The identities in (\ref{e19}) follow from
\begin{eqnarray} \label{e21}
(R_n^{|m|})'(r) & = & \sum_{j=0}^{\frac12(n-1-|m-1|)}\,(n-2j)\,R_{n-1-2j}^{|m-1|}(r) \nonumber \\[3.5mm]
& & +~\sum_{j=0}^{\frac12(n-1-|m+1|)}\,(n-2j)\,R_{n-1-2j}^{|m+1|}(r)~,
\end{eqnarray}
see \cite{ref4}, (10) and add the $\pm$-cases that occur there, where we also recall that $R_{n-1-2j}^{|m\pm 1|}(1)=1$ for $j$ in the summation ranges in the two series in (\ref{e21}). Actually, (\ref{e21}) was used in \cite{ref6}, Appendix~B, to show the bound in (\ref{e2}) by induction. In general, it can be shown by induction from (\ref{e21}) that the ratio of $(R_n^{|m|})^{(k)}(1)$ and the right-hand side of (\ref{e3}) converges to $(1/2)_k/(1)_k$ as $n\pr\infty$ and $m$, $k$ are fixed.

In the case that $k=n-2$, we have
\begin{eqnarray} \label{e22}
(R_n^{|m|})^{(n-2)}(1) & = & a_n(m)\,T_n^{(n-2)}(1)+a_{n-2}(m)\,T_{n-2}^{(n-2)}(1) \nonumber \\[3mm]
& = & a_n(m)\,T_n^{(n-2)}(1)\Bigl(1+\frac{|m|^2}{n}~\frac{1}{n(2n-3)}\Bigr)~,
\end{eqnarray}
where the first item in (\ref{e10}) and the second form in (\ref{e18}) for $T_{n-2}^{(n-2)}(1)$ and $T_n^{(n-2)}(1)$ has been used. Hence, the lower bound in (\ref{e15}) is also sharp for $k=n-2$ when $n\pr\infty$, and a similar thing can be concluded for $k=n-4$ from the second item in (\ref{e10}).

We next consider $a_n(m)$ as given by (\ref{e9}) for $n>0$. We have for $0\leq m\leq\frac12\,n$ by Stirling's formula
\beq \label{e23}
a_n(m)=\frac{4}{\sqrt{2\pi n}}\,\Bigl(\frac{1}{1-\dfrac{|m|^2}{n^2}}\Bigr)^{\frac{n+1}{2}}\,
\Bigl(\frac{1-\dfrac{|m|}{n}}{1+\dfrac{|m|}{n}}\Bigr)^{\frac{|m|}{2}}\, \Bigl(1+O(\dfrac1n\Bigr)\Bigr)~.
\eq
In the case that $|m|=O(\sqrt{n})$, this gives
\beq \label{e24}
a_n(m)=\frac{4}{\sqrt{2\pi n}}\,\exp\Bigl({-}\,\frac{|m|^2}{2n}\Bigr)\,\Bigl(1+O(\frac{1}{\sqrt{n}}\Bigr)\Bigr)~.
\eq
Therefore, by (\ref{e15}), the upper bound (\ref{e3}) is attained within a factor $O(1/\sqrt{n})$ when $|m|=O(\sqrt{n})$.

There are also cases that neither bound in (\ref{e15}) is sharp. For instance, when $n=m=2k$ with large $k$, the ratio of $(R_n^{|m|})^{(k)}(1)$ and the lower bound and upper bound in (\ref{e15}) are of the order $(32/37)^k$ and $(8/27)^k$, respectively, and thus tend to $\infty$ and 0 exponentially fast as $k\pr\infty$. \\ \\
{\bf Acknowledgement}. \\
The author wishes to thank B.M.M.\ de Weger for demonstrating non-triviality of the problem of deciding whether $P_p^{(\gamma,\delta)}(0)=0$ (non-negative integer $\gamma$, $\delta$ and $p$), already for degrees $p$ as low as 4, 5, 6. \\ \\ \\ \\
\mbox{}

\end{document}